\documentclass[
	a4paper,
	11pt,
	headsepline,
	numbers=endperiod,
	]{scrartcl}
\usepackage{geometry}
\usepackage{caption}
\setkomafont{sectioning}{\normalfont\normalcolor\bfseries}
\usepackage{mathrsfs}
\usepackage{braket}
\usepackage{setspace}
\spacing{1.2}
\usepackage{microtype}
\usepackage{graphicx}
\usepackage{amsmath,amsfonts,amssymb,xfrac,nicefrac}
\usepackage{aliascnt}					
\usepackage{amsthm}
\usepackage{mathtools}
\usepackage{hyperref}
\usepackage[capitalise]{cleveref}

\theoremstyle{definition}
\newtheorem{correspondence}{Correspondence}[section]
\newtheorem{definition}[correspondence]{Definition}
\newtheorem{example}[correspondence]{Example}

\begin{document}

\title{\LARGE{Dimension Reduction of Two-Dimensional Persistence via Distance Deformations}}

\author{Maximilian Neumann\footnote{Email: neumann.mn@icloud.com}}
\date{ }
\maketitle

\noindent \textbf{Abstract.} This article grew out of the application part of my Master’s thesis at the Faculty of Mathematics and Information Science at Ruprecht-Karls-Universität Heidelberg under the supervision of PD Dr.~Andreas Ott. In the context of time series analyses of RNA virus datasets with persistent homology, this article introduces a new method for reducing two-dimensional persistence to one-dimensional persistence by transforming time information into distances.

\addsec{Introduction}

Let $(S,h)$ be a finite distance space where $h$ takes values in $\mathbb{N}$. 
Assume that we have a time dependent filtration
\begin{equation}
S_0 \subseteq \dots \subseteq S_m=S \notag
\end{equation}
where $m$ is the number of time steps (e.g.~months or days). For $i \in \{0, \dots, m\}$, consider the Vietoris-Rips filtration
\begin{equation}
V_i:=(V_{0}(S_i,h) \subseteq V_{1}(S_i,h) \subseteq V_{2}(S_i,h) \subseteq \dots ). \notag
\end{equation}
As in the work of Bleher et al.~\cite{topologyidentifies2021}, where $S$ is a finite set of SARS-CoV-2 RNA sequences with Hamming distance $h$, we are interested in detecting cycles that correspond to bars in the barcode $\mathcal{B}(H_1(V_i))$ born in the first filtration step. In \cite{topologyidentifies2021}, these cycles are called \textit{single nucleotide variation} (SNV) cycles and are used for a topological recurrence (time series) analysis of SARS-CoV-2, where each time step $\mathcal{B}(H_1(V_i))$ is computed with Ripser \cite{bauer2021ripser}. Ripser is a highly optimised software tool for the computation of persistent homology, capable of processing hundreds of thousands of distinct RNA sequences \cite{topologyidentifies2021}. However, this classical approach, where each time step is computed seperately, can become very time-consuming for large $m$. In this article, we present a new method that improves this classical approach and enables the extraction of SNV cycles for each time step with only one barcode computation.

The $V_i$ naturally lead to a finite bifiltered simplicial complex $V$. The f.g.~two-dimensional persistence module $H_1(V)$ contains all the information that occur within the $H_1(V_i)$. Moreover, $H_1(V)$ contains additional information about the behaviour of homology classes along the time filtration parameter. In order to track these time-persistence features, we introduce a distance deformation technique to compute the barcode of a certain one-dimensional subfiltration $R\subseteq V$, which is relevant for the detection of SNV cycles. For this, we deform $h$ into a new distance $h^*$ on $S$ to realise $R$ as a Vietoris-Rips filtration $R^*$, such that we have a correspondence between the barcodes $\mathcal{B}(H_1(R^*))$ and $\mathcal{B}(H_1(R))$ for the bars corresponding to SNV cycles. Since $R^*$ is a Vietoris-Rips filtration, the barcode $\mathcal{B}(H_1(R^*))$ can be computed with Ripser \cite{bauer2021ripser}. In practical experiments, one could investigate whether this new method provides a performance advantage over the classical approach to a time series analysis, where each time step is computed seperately.
\vspace{0.5cm}

\noindent \textbf{Acknowledgements.} I would like to take this opportunity to thank Andreas Ott and Michael Bleher for their support and the inspiring discussions.

\section{SNV cycles}
For this article, let $(S,h)$ be a finite distance space, i.e. $S$ is a finite set and $h$ is a metric or more generally a semimetric\footnote{A semimetric satisfies all the axioms of a metric with exception of the triangle inequality.} on $S$. 
Recall that for every $r \in \mathbb{R}$, the \textit{Vietoris-Rips complex of} $(S,h)$ at \textit{scale} $r \in \mathbb{R}_{\geq 0}$ is the abstract simplicial complex
\begin{equation}
V_{r}(S,h):=\{ \emptyset \neq \sigma \subseteq S \mid  h(x,y)\leq  r \, \forall x,y \in \sigma \} \notag.
\end{equation}
Assume that we have a (time dependent) filtration
\begin{equation}
S_0 \subseteq \dots \subseteq S_m=S \notag. 
\end{equation}
For $i \in \{0, \dots, m\}$, consider the Vietoris-Rips filtration
\begin{equation}
V_i:=(V_{0}(S_i,h) \subseteq V_{1}(S_i,h) \subseteq V_{2}(S_i,h) \subseteq \dots ). \notag
\end{equation}
Denote by $H_1(V_i)=H_1(V_i, \mathbb{F}_p)$ the first simplicial homology with coefficients in a finite prime field $\mathbb{F}_p$ applied to the filtration $V_i$. Then $H_1(V_i)$ is a finitely generated (f.g) one-dimensional persistence module.

As in the work of Bleher et al.~\cite{topologyidentifies2021} where $S$ is a finite set of SARS-CoV-2 RNA sequences with Hamming distance $h$, we are interested in detecting cycles that correspond to bars in the barcode $\mathcal{B}(H_1(V_i))$ born in the first filtration step. In \cite{topologyidentifies2021}, these cycles are called \textit{single nucleotide variation} (SNV) cycles and are used for a topological recurrence (time series) analysis of SARS-CoV-2. For simplicity, we also call such cycles SNV cycles in our more general setting. 
\begin{definition}[SNV cycle] The underlying homology class representatives of bars in the barcode $\mathcal{B}(H_1(V_i))$ born in the first filtration step $H_1(V_1(S_i,h))$ are called \textit{SNV cycles} in time step $i$. 
\end{definition}
For every $i \in \{0, \dots, m\}$, denote by $\mathbf{SNV}_i$ a full set of SNV cycle representatives extracted from the barcode $\mathcal{B}(H_1(V_i))$. In \cite{topologyidentifies2021}, the barcodes $\mathcal{B}(H_1(V_i))$ are computed with Ripser \cite{bauer2021ripser} and the $\mathbf{SNV}_i$ are extracted from the Ripser output. Ripser is a highly optimised software tool, capable of processing hundreds of thousands of distinct RNA sequences \cite{topologyidentifies2021}. However, this classical approach to a time series analysis has the following issues:

\begin{enumerate}
\item Computing each time step seperately can be very time consuming for large $m$ (e.g. a time series analysis over one year on a daily basis).
\item We are not able to track the time-stability of SNV cycles, i.e. whether the image of the homology class $[\omega]$ of an SNV cycle $\omega \in \mathbf{SNV}_i$ under the canonical homomorphism
\begin{equation}
H_1(V_1(S_i,h)) \longrightarrow H_1(V_1(S_{i+1},h)) \notag
\end{equation}
is zero or not.

\item Since each time step is computed  seperately, the $\mathbf{SNV}_{i}$ are not automatically compatible: let $\omega \in \mathbf{SNV}_i$ and assume that the image of $[\omega]$ under the canonical homomorphism
\begin{align}
H_1(V_1(S_i,h)) \longrightarrow H_1(V_1(S_{i+1},h)) \notag
\end{align}
is not zero. Then it still may happen that $\omega \not \in \mathbf{SNV}_{i+1}$.

\end{enumerate}
In Sections \ref{sec:Multipersistence dimension reduction for SNV cycles} and \ref{sec:Hamming distance deformation}, we present a method that enables the extraction of SNV cycles for each time step with only one barcode computation. The resulting SNV cycles are automatically compatible and we can track their time-stability.
\section{Dimension reduction}\label{sec:Multipersistence dimension reduction for SNV cycles}

The $V_i$ naturally lead to a finite bifiltered simplicial complex $V$. We obtain a f.g.~two-dimensional persistence module $H_1(V)$ which contains all the information that occur within the $H_1(V_i)$. Moreover, $H_1(V)$ contains additional information about the behaviour of homology classes along the time filtration parameter. Since we are only interested in detecting SNV cycles and not in determining their lifespan in the barcodes $\mathcal{B}(H_1(V_i))$, it suffices to compute the barcode $\mathcal{B}(H_1(R))$  where  $R\subseteq V$  is the one-dimensional subfiltration
\begin{align}
R_i :=  \begin{cases} V_{0}(S_i,h), & i=-1 \notag \\
V_{1}(S_i,h), & i \in  \{0, \dots, m\} \\
V_{i-m+1}(S_m,h),  & i \in \mathbb{N}_{\geq m+1} \notag
\end{cases}
\end{align}
For reasons of notation, we start with $i=-1$. The f.g.~one-dimensional persistence module $H_1(R)$ can be viewed as a dimensional reduction of $H_1(V)$. The barcode $\mathcal{B}(H_1(R))$ contains all the information we need to extract SNV cycles for each time step $i\in  \{0,  \dots, m\}$. Moreover, $\mathcal{B}(H_1(R))$ tracks the stability of SNV cycles along the time filtration parameter.  

The idea to consider barcodes of subfiltrations follows a more general concept introduced by Carrie et al.~\cite{Bettinumbersmultipers} and called fibered barcode by Lesnick and Wright \cite{lesnick2015interactive}. Fibered barcodes are closely related to the rank invariant introduced by Carlsson and Zomorodian in \cite{Carlsson2009multidimensionalpersistence}. In \cite{Bettinumbersmultipers}, it is shown that the fibered barcode and the rank invariant determine each other.

\section{Distance deformation}\label{sec:Hamming distance deformation}
In this section, we introduce a distance deformation technique to realise $R$ as a Vietoris-Rips filtration $R^*$ such that we have a correspondence between the barcodes $\mathcal{B}(H_1(R^*))$ and $\mathcal{B}(H_1(R))$ for the bars corresponding to SNV cycles.

For the following, let $N=N(m)$ be the lowest power of $10$ such that $m < N$. For example, if $m=34$, then $N=100$. For $x \in S$, let \begin{equation}
D(x) := \min \{ i \in \{0, \dots, m\} \mid x \in S_i \} \notag.
\end{equation}

\begin{definition}[Distance deformation] We define a new distance $h^*$ on $S$ as follows: let $x,y \in S$ with $D(x) \geq D(y)$. Define
\begin{align}
h^*(x,y) := \begin{cases}
h(x,y)+ D(x)/N, & x \neq y  \\
0, & x=y \end{cases}
\notag
\end{align}
and \begin{equation}
h^*(y,x):= h^*(x,y) \notag.
\end{equation}
\end{definition}

\begin{example}The intuition behind $h^*$ is that time information is transformed into distances. Let $m=364$. Then we have $N=N(364)=1,000$. Let $x,y,z \in S$ with $D(x)=264$ and $D(y)=D(z)=132$. Assume that $h(x,y)=h(x,z)= h(y,z)=1$.
Then we have \begin{equation}h^*(x,y)=h^*(x,z)=1.264 \notag
\end{equation}
and 
\begin{equation}
h^*(y,z)=1.132 \notag.
\end{equation}
\end{example}
\begin{figure}
\captionsetup{format=plain, labelfont={bf},labelformat={default},labelsep=none,name={Figure}}
\centering
\includegraphics[scale=0.4]{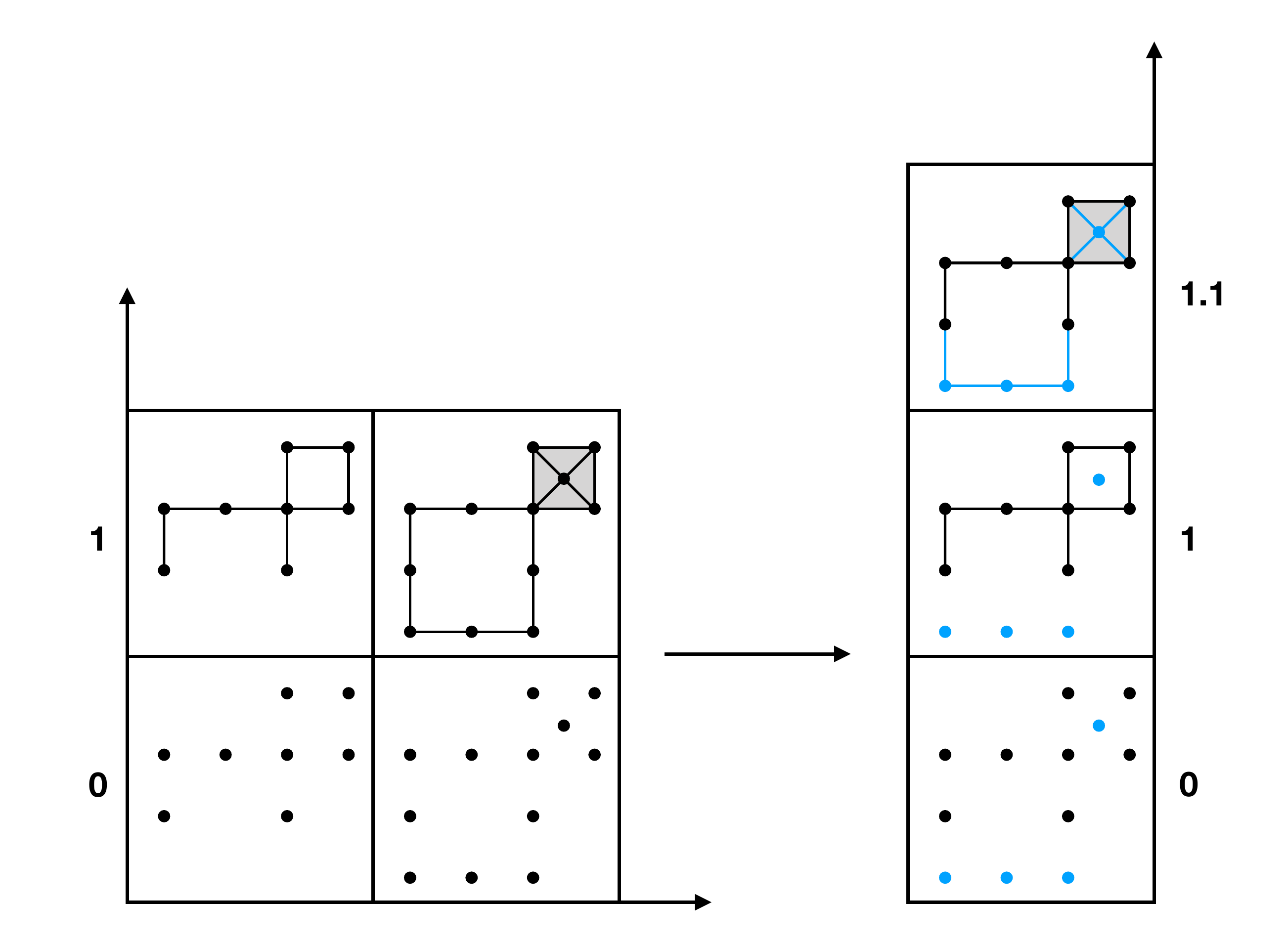}
\caption{\, Here we illustrate the correspondence between SNV cycles and their deformed equivalents. The blue-coloured points and edges indicate that the distance was deformed according to the time step they were added. As we can see, an SNV cycle was destroyed by adding a point along the time filtration parameter.}
\label{fig:metricdeformation}
\end{figure}
Consider the Vietoris-Rips filtration $R^*$, where for $i \in \mathbb{Z}_{\geq -1}$,
\begin{align}
R^*_i:=V_{\kappa_i }(S,h^*)  \notag
\end{align}
with filtration parameters
\begin{equation}
\kappa_i:= \begin{cases} 0, & i=-1  \\
1+ i/N, & i \in\{0, \dots, m\} \\
2+ (i-(m+1))/N, & i \in \{m+1, \dots, 2m+1 \}\\
\dots
\end{cases} \notag
\end{equation}
Then $H_1(R^*)$ is a f.g.~one-dimensional persistence module. By construction, we have the following correspondence (illustrated in Figure \ref{fig:metricdeformation}).

\begin{correspondence}\label{correspondence} Consider the barcodes $\mathcal{B}(H_1(R))$ and $\mathcal{B}(H_1(R^*))$. Let $i \in \{0, \dots, m\}$. Then bars born in $H_1(R^*_i)$ are in one to one correspondence with bars born in $H_1(R_i)$. Let $j \in \{0, \dots, m-i\}$. If a bar born in $H_1(R_i)$ dies in $H_1(R_{i+j})$, the corresponding bar born in $H_1(R^*_i)$ dies in $H_1(R^*_{i+j})$. 
\end{correspondence}
Using this correspondence, the definition of SNV cycles translates as follows.
\begin{definition}[Deformed SNV cycle] The underlying homology class representatives of bars in the barcode $\mathcal{B}(H_1(R^*))$ born in $H_1(R^*_0), \dots, H_1(R^*_m)$ are called \textit{deformed SNV cycles}.
\end{definition}
Denote by $\mathbf{SNV}^*$ a full set of deformed SNV cycle representatives extracted from $\mathcal{B}(H_1(R^*))$. For ${i \in \{0, \dots, m\}}$, define \begin{equation}\mathbf{SNV}^*_{i} := \{ \omega \in  \mathbf{SNV}^* \mid 0 \neq [\omega]  \in H_1(R^*_{i}) \}   \notag.
\end{equation}
By construction, we have a bijection of sets
\begin{equation}
\mathbf{SNV}^*_i \cong \mathbf{SNV}_i \notag.
\end{equation}
Moreover, we have compatibility: let $\omega \in \mathbf{SNV}^*_i$ and assume that the image of $[\omega]$ under the canonical homomorphism
\begin{equation}
H_1(R^*_i) \longrightarrow H_1(R^*_{i+1}) \notag
\end{equation}
is not zero. Then $\omega \in \mathbf{SNV}^*_{i+1}$ by construction. In addition, we can track the time-stability of SNV cycles and instead of $m$ barcode computations of $\mathcal{B}(H_1(V_i))$ for $i \in \{0, \dots,m\}$, only the computation of $\mathcal{B}(H_1(R^*))$ has to be performed. Since $R^*$ is a Vietoris-Rips filtration, the barcode $\mathcal{B}(H_1(R^*))$ can be computed with Ripser \cite{bauer2021ripser}. In practical experiments, one could investigate whether this new method provides a performance advantage over the classical approach to a time series analysis, where each time step is computed seperately.

\bibliographystyle{plain}
\bibliography{article}

\end{document}